\documentclass[12pt]{amsart}
\usepackage{fullpage, amsmath, amsthm,amsfonts,amssymb,stmaryrd, mathrsfs}
\usepackage{hyperref}

\usepackage[pdftex]{color}
\usepackage[all]{xy}
\usepackage{caption}
\usepackage{graphicx}

\numberwithin{equation}{subsection}
\newtheorem{theorem}[subsection]{Theorem}

\renewcommand{\qed}{$\blacksquare$}
\newcommand{\dia}{$\diamond$}

\newtheorem{lemma}[subsection]{Lemma}

\theoremstyle{definition}

\newtheorem{ex}[subsection]{Example}
\newtheorem{rem}[subsection]{Remark}

\def\FF{\mathbb{F}}

\def\PP{\mathbb{P}}
\def\QQ{\mathbb{Q}}

\def\ZZ{\mathbb{Z}}

\begin{document}

\title{Computing zeta functions of large polynomial systems over finite fields}
\author{Qi Cheng}
\address{University of Oklahoma, School of
Computer Science, Norman, OK 73019}
\email{qcheng@ou.edu}
\author{J. Maurice Rojas}
\address{Texas A\&M University, Department of
Mathematics, College Station, TX 77843-3368}
\email{rojas@math.tamu.edu}
\author{Daqing Wan}
\address{University of California, Irvine, Department of
Mathematics, Irvine, CA 92697-3875}
\email{dwan@math.uci.edu}

\date{\today}

\begin{abstract} In this paper, we improve the algorithms of Lauder-Wan 
\cite{LW} and Harvey \cite{Ha} to compute the zeta function of a system 
of $m$ polynomial equations in $n$ variables over the finite field $\FF_q$ of 
$q$ elements, for $m$ large. The dependence on $m$ in the original algorithms 
was exponential in $m$.  Our main result is a reduction of the exponential 
dependence on $m$ to a polynomial dependence on $m$. As an application, we 
speed up a doubly exponential time algorithm from a software verification 
paper \cite{BJK} (on universal equivalence of programs over finite 
fields) to singly exponential time. One key new ingredient is 
an effective version of the classical Kronecker theorem which 
(set-theoretically) reduces the number of defining equations for a ``large''  
polynomial system over $\FF_q$ when $q$ is suitably large. 
\end{abstract}

\maketitle

\setcounter{tocdepth}{1}

\section{Introduction}

Let $\FF_q$ be the finite field of $q$ elements with characteristic $p$. Let 
$F$ be a polynomial system with $m$ equations and $n$ variables over 
$\FF_q$: 
$$F(x_1,\ldots, x_n)=(f_1(x_1,\ldots, x_n),\ldots, f_m(x_1,..., x_n)),$$
where each $f_i \in \FF_q[x_1,\ldots, x_n]$ is a polynomial in $n$ variables of degree at most $d$. Note that the total number of digits 
needed to write down the monomial term expansion of such a system is 
$O(m( d+1 )^n\log q)$. So it is 
natural to use $m( d +1 )^n\log q$ as a measure of {\em input size} for $F$ when 
we start discussing algorithmic efficiency.  
For our purposes here, and for reasons to be made clear shortly, we will 
call the polynomial system $F$ {\em large} if the number $m$ of equations is 
at least $n+2$.  

A basic algorithmic problem in number 
theory is to compute the number $N_q(F)$ of solutions of the polynomial system $F=(0,\ldots,0)$ over $\FF_q$. More precisely, we set 
$$N_q(F) := \left\{\left. (x_1,\ldots, x_n)\in \FF_q^n \; \right| \; 
F(x_1, \ldots, x_n)=(0,\ldots,0)\right\}.$$ 
The special case $(m,n)\!=\!(1,2)$ already plays a huge role in cryptography, 
since curves with a specified number of points are crucial to the design of 
many cryptosystems (see, e.g., \cite{handbook}). 
 
An even deeper basic problem is to consider all extension fields of 
$\FF_q$ and compute the full sequence $N_q(F), N_{q^2}(F),\ldots, N_{q^k}(F), \ldots$ or, equivalently, the generating zeta function 
$$Z(F, T) =\exp\!\left(\sum_{k=1}^{\infty} \frac{N_{q^k}(F)}{k}T^k\right).$$ 
Understanding this generating function occupied a good portion of 20th 
century algebraic and arithmetic geometry. Interestingly, this generating 
function has found a recent application to software engineering, specifically, 
in program equivalence \cite{BJK}. (We clarify this in the next section.) 
It is not at all obvious from the definition that this zeta function is 
effectively computable, so let us briefly recall how it actually is. 

A deep and celebrated theorem of Dwork from 1960  says that the zeta function 
is a rational function in $T$. A theorem of Bombieri \cite{Bo} from 1988 says 
that the total degree of the zeta function is effectively bounded. 
It then follows, from basical manipulation of power series, that the zeta 
function is effectively computable, although practical efficiency is far 
more subtle: See \cite{Wa} for a survey on algorithms for computing zeta 
functions. A general deterministic algorithm to compute $Z(F,T)$ was 
constructed in Lauder-Wan \cite{LW} with running time 
$$2^m (pm^nd^n\log q)^{O(n)}.$$ 
For small characteristic $p$, this general purpose algorithm remains the best 
so far. However, for large characteristic $p$, 
the dependence on $p$ has been improved by Harvey \cite{Ha}, who constructed 
an algorithm with running time 
$$2^m p (m^nd^n\log q)^{O(n)}.$$ 
(There is also a variant in \cite{Ha} with time complexity linear in 
$\sqrt{p}$ instead, but at the expense of increasing the space complexity to 
roughly the same order as the time complexity.) 
The algorithms from \cite{LW} and \cite{Ha} are, however, fully exponential in 
$m$, even for fixed $n$. 

To improve the dependence on $m$,  
we briefly explain how the exponential factor $2^m$ arises in the 
algorithms of \cite{LW} and \cite{Ha}. Both algorithms, in the case $m=1$ 
(the hypersurface case), are obtained via 
$p$-adic trace formulas (meaning linear algebra with large matrices over 
the polynomial ring $(\ZZ/p^\lambda \ZZ)[t]$, arising after some cohomological 
calculations). The case $m>1$ is then reduced to the case $m=1$ via an 
inclusion-exclusion trick \cite{Wa} to compute the 
zeta function for each of the $2^m$ hypersurfaces defined by 
$f_S =\prod_{i\in S} f_i$, 
where $S$ runs through all subsets of $\{1,2,\ldots, m\}$ and $\deg(f_S) \leq 
|S|d\leq md$. 

In this paper, we improve the 
Lauder-Wan algorithm and the Harvey algorithm by using a different reduction to reduce the exponential factor $2^m$ to $m$. 
One key new idea is to prove an effective version of Kronecker's theorem which reduces the number 
$m$ of defining equations to $n+1$ if $q$ is suitably large: See 
Section \ref{sec:kronecker} below.  

Our main result is the following: 
\begin{theorem} 
\label{thm:main} 
There is an explicit deterministic algorithm which computes the zeta function 
$Z(F,T)$ of the system $F$ over $\FF_q$ (with $m$ equations, $n$ variables,  
of degree at most $d$) in time $$mp(n^nd^n\log q)^{O(n)}.$$
\end{theorem}

We will see in the next section how our theorem enables us to speed up a 
{\em doubly} exponential time algorithm (from \cite{BJK}) for program 
equivalence to {\em singly} exponential time. In particular, we will now 
briefly review some of the background on programs over finite fields.  

\section{Programs, Their Equivalence, and Zeta Functions} 
A basic and difficult problem from the theory of programming languages 
is determining when two programs always yield the same output (hopefully 
without trying all possible inputs). This 
problem --- a special case of {\em program equivalence} --- also has an 
obvious parallel in cryptography: a fundamental problem  
is to decide whether a putative key for an unknown stream cipher (that one has 
spent much time decrypting) is correct or not, without trying all possible 
inputs. In full generality, program equivalence is known to be undecidable in 
the classical Turing model of computation. However, program equivalence (and 
{\em formal verification}, in greater generality \cite{dagstuhl}) 
remains an important need in software engineering and cryptography. 
It is then natural to ask these questions in a more limited setting.  

For instance, Barthe, Jacomme, and Kremer (in \cite{BJK}) describe 
a programming language which enables a broad family of calculations 
(and verifications thereof) involving polynomials over finite fields. 
They proved that program equivalence in their setting is decidable, and 
gave an algorithm with doubly exponential complexity. We now briefly 
review their terminology (from \cite[Sec.\ 2.2]{BJK}), and how their algorithm 
requires a non-trivial use of zeta functions.  

To be more precise, in their restricted setting, a {\em program} is a 
sequence of logical/polynomial expressions over a finite field. To 
define this rigorously, one first 
fixes a set $I$ of {\em input} variables 
and a set $R$ of {\em random} variables.
Then all 
possible expressions making up a program can be defined recursively 
(building up from (1) and (2) below) as follows:  
\begin{enumerate}
\item a polynomial $P\in\FF_{q}[I,R]$; 
\item the failure statement $ \perp $;
\item an ``{\tt if}'' statement of the following form:
\[ \texttt{\ if\ } b \texttt{\ then\ } e_1 \texttt{\ else\ } e_2  \]
where $ e_1 $ and $ e_2 $ are expressions, and $ b $ is a propositional 
logic formula, whose atoms are of the form $ Q=0 $ for some 
$ Q\in\FF_{q}[I,R]$.
\end{enumerate}

\begin{rem}{\em  
  Programs in \cite{BJK} are written using semi-colons as delimiters, 
  similar to some real-world program languages
such as C or Java.\dia } 
\end{rem} 

The {\em size} of a program is defined to be the number of characters in a 
program. The presence of random variables enables our 
programs to use randomization, and give answers with a certain probability 
of failure. We denote the set of all such programs by 
$ {\mathcal{P}_{q}} (I,R) $.  Polynomials in a
program are represented by arithmetic formulas, so the degree of any polynomial
in the program is bounded from above by the size of the program.
Note that programs in this core
language do not have loops. If a program has neither ``{\tt if}'' statements 
nor failure statements then we call the program an {\em arithmetic program}. 
The set of all arithmetic programs is denoted by 
$ \bar{\mathcal{P}}_{q} (I,R) $.

The {\em number of expressions at the top level} of a program $ \PP $ --- 
denoted by $ |\PP| $ --- is simply the length of the sequence 
defining $\PP$. (In a real-world programming language, the ``top level'' 
of a program simply means one ignores subroutines and, e.g., statements 
{\em inside} of an ``if'' statement.) Note also that since our programs can 
use random variables, our programs thus send input values 
in $ \FF_{q^k}^{|I|} $ to a probability distribution over 
$ \FF_{q^k}^{|\PP|} $ for any positive integer $ k $.  
Assuming that the program does not fail (i.e., there is no evaluation of  $\perp$ that halts the program),
this can be viewed as the following map of inputs to maps:
$\FF_{q^k}^{|I|} \rightarrow (   \FF_{q^k}^{|\PP|} \rightarrow [0,1]  ) $.  
It is clear that understanding the semantics of a program requires counting
solutions of a polynomial system.

\begin{ex} {\em 
Fixing $I= \emptyset $ and  $  R\!=\!\{x\}$, the program \\ 
\mbox{}\hfill $x*x$ ; $x*x*x$ \hfill \mbox{} \\ 
outputs a uniformly random square and a uniformly random cube from 
$\FF_{q}$, though these two numbers are not independent.\footnote{We use $x*x$ in place of $x^2$ since
  polynomials are represented by arithmetic formulas.
}
Let $ N(\alpha,\beta) $ denote  the number of solutions in $ \FF_{q^k} $ of
\begin{align*}
  x^2 &= \alpha\\
  x^3 &= \beta
\end{align*}
The program outputs a distribution sending $ (\alpha,\beta)\in \FF_{q^k}^2 $ to
$ N(\alpha,\beta)/q^k $.
\dia}
\end{ex} 

\begin{ex} {\em 
Let $ I=\{ x\} $ and $ R = \{y,z\} $. The following  program 
$ \PP_1 $ is in $ \mathcal{P}(I,R) $    
\[  \mathtt{\ if\ } \neg (x =0) \mathtt{\ then\ } y+1 \mathtt{\ else\ } y+2;\ \
     z*z  \]
The program $\PP_1$ yields the probability distribution on 
$\FF^2_{q^k}$ corresponding to the first coordinate being uniformly random in 
$\FF_{q^k}$ and the second coordinate a uniformly random square in 
$\FF_{q^k}$. \dia}
\end{ex}
To calculate the distribution, the sample space consists of the assignments
to random variables so that the program does not fail. For example,  
the following program ($ I=\{x\} $ and $ R=\{y\} $  ) computes the inverse of $ x $  with probability $ 1 $:
\[ \texttt{\ if\ }  x=0 \texttt{\ then\ } 0 \texttt{\ else\ } \texttt{\ if\ }
  x*y =1 \texttt{\ then\ } y \texttt{\ else\ } \perp \]

Given two programs, we would like to check whether they produce the same 
distribution for any input. More generally,
let $ \PP_1, \QQ_1 $ be programs and $ \PP_2, \QQ_2 $ be  arithmetic programs. 
We write $ \PP_1 | \PP_2 \approx \QQ_1 | \QQ_2 $ if,    
 taking any input $ c $ under the condition that $ \PP_2 = \vec{0}$,
$ \PP_1 $ outputs the same distribution as  $ \QQ_1 $ taking $ c $ as input 
under the condition $ \QQ_2 = \vec{0}$. To calculate the distribution,
we only need to consider the random values such that
none of $ \PP_1 $ and $ \QQ_1 $ output $ \perp $.   
\begin{rem} \label{rem:zeta} 
{\em  
Observe that the set of inputs yielding a fixed sequence of outputs is 
nothing more than a constructible set over $\FF_{q^k}$, i.e., a boolean 
combination of algebraic sets over $\FF_{q^k}$. In particular, the 
set of inputs making two programs {\em differ} is also a constructible 
set over a finite field. \dia}  
\end{rem}  

The question of equivalence can be raised for a fixed $ k $, or
for all positive integers $ k $. The latter case is called {\em universal
equivalence}, which is most relevant to our discussion here. For example, 
let $ \QQ_1 $ be the program defined by: 
\[ y;  \ \ 7*(z+1)*(z+1).  \]
If $ 7 $ is a nonzero square in $ \FF_{q} $ then $ \PP_1 $ is universally
equivalent
to $ \QQ_1 $, i.e., $ \PP_1|0 \approx \QQ_1|0 $.
Otherwise, $\PP_1|0$ and $\QQ_1|0$ are not equivalent over $ \FF_{q} $,
and hence not universally equivalent.

\subsection{How \cite{BJK} reduces from general to arithmetic programs} 
Note that in greater generality, checking universal equivalence means 
checking if a {\em sequence} of constructible sets consists solely 
of empty sets (per Remark \ref{rem:zeta} above). As observed in 
\cite{BJK}, this can be done by a single zeta function computation. 
This is, in essence, how \cite{BJK} proved that universal equivalence for 
arithmetic programs can be done in singly exponential time. 
For universal equivalence of conditional programs, the same ideas apply, but 
\cite{BJK} proved a doubly exponential complexity upper bound. More precisely, 
for general programs $ \PP_1, \QQ_1 $ and arithmetic programs $ \PP_2, \QQ_2$ ,
they defined a reduction,  
to obtain four arithmetic programs $ \PP_1', \PP_2', \QQ_1' $ and $ \QQ_2' $ 
so that
\[ \PP_1|\PP_2 \approx \QQ_1 | \QQ_2 \mathbf{\ if\ and\ only\ if\ } \PP_1'|\PP_2' \approx \QQ_1'
  | \QQ_2' \]

It is clear that one needs to be able to remove failure statements 
($  \perp  $) and ``{\tt if}'' statements from $ \PP_1 $ ( and repeat
the procedures on $ \QQ_1 $ )
in order for such a reduction to work.
Here, we will use examples to illustrate the ideas in the reduction.
See \cite{BJK} for the full, formal treatment.

We may assume that there is at most one occurrence of the failure statement
in $ \PP_1 $, since we can collect the conditions for
failure together. For example, the following program
\[ \texttt{\ if\ }  A_1 \texttt{\ then\ } \perp \texttt{\ else\ } P_1;
\texttt{\ if\ }  A_2 \texttt{\ then\ } P_2 \texttt{\ else\ } \texttt{\ if\ }  A_2 \texttt{\ then\ } \perp \texttt{\ else\ } P_3 \]
is equivalent to 
\[ \texttt{\ if\ }  A_1 \vee (\neg A_2 \wedge A_2 ) \texttt{\ then\ } \perp \texttt{\ else\ } P_1;
\texttt{\ if\ }  A_2 \texttt{\ then\ } P_2 \texttt{\ else\ }  P_3 \]
The new program has length polynomial in the length of the old program, since 
the number of $ \perp $ in the input
program is bounded from above by the length of the input.
Without loss of generality, suppose that $ \PP_1 $ has the form  
{\par\centering \tt if b then $ P_1 $  else $ \perp $; $\cdots$ \par}
\noindent where $ \perp $ occurs only once in the program. 
If the condition $ b $ is a disjunction of literals\footnote{A 
{\em disjunction} is simply a boolean ``OR'' applied to several 
propositions. A {\em literal} 
is simply a variable, or the negation thereof.} then we can find a single
polynomial $ B $ whose vanishing represents $ b $. For example, if $ b $ is 
\mbox{$ (P_2=0)\vee \neg (P_3=0) \vee \neg (P_4=0)$,} then 
we construct the polynomial $B=$\mbox{$P_2 (t_3 P_3-1) (t_4 P_4-1  ) $.}
The new programs become 
\begin{align*}
  \PP'_1 &= P_1; \cdots \\
 \PP'_2 &=   \PP_2;  B; t_3 (t_3 P_3-1);
  P_3 (t_3 P_3-1); t_4 (t_4 P_4-1); P_4 (t_4 P_4-1)
 \\ 
  \QQ'_1 &=\QQ_1\\
  \QQ'_2 &= \QQ_2; B; t_3 (t_3 P_3-1);
  P_3 (t_3 P_3-1); t_4 (t_4 P_4-1); P_4 (t_4 P_4-1)
\end{align*}
Here $ t_3 $ and $ t_4 $ are new random variables but they are uniquely
determined by $ P_3 $ and $ P_4 $ under the constraints.  
Namely if $ P_3=0 $, then $ t_3 =0 $, otherwise $ t_3 = 1/P_3 $.   
For a more general 
proposition formula $ b $, we first convert it  to a CNF formula,\footnote{{\em 
Conjunctive Normal Form}, meaning ``an AND of ORs''...} which may
result in a conjunction of
exponentially many disjunctions, hence exponentially many polynomials $ B_1, B_2,
\cdots, B_m $, in addition to polynomials like $ t_i (t_i P_i-1) $
and $ P_i(t_i P_i-1) $ etc. The new equivalence is
\[  P_1; \cdots | \PP_2, B_1, B_2,
\ldots \approx \QQ_1 | \QQ_2, B_1, B_2,
\ldots   \]
Nevertheless we only introduce polynomially many new variables,
since we need at most one new variable for each
polynomial in the original program. Also the $ \PP_2' $ may be exponentially
long, but the $ \PP_1' $ is actually shorter than the original $ \PP_1 $.\\ 

Observe that we may also assume that
all the inputs to conditional statements are literals. 
For example we can replace
{\par\centering\tt if $A_1  \vee  A_2\ $ then $P_1 $ else $P_2$ \par }
by
{\par\centering\tt if $A_1$ then $P_1$ else if $A_2$ then $ P_1 $ else $ P_2 $. \par }
Then, to remove ``{\tt if}'' in a conditional
statement such as  
\[ \cdots; \texttt{\ if\ }  \neg ( B=0 )\texttt{\ then\ } P_1 \texttt{\ else\ } P_2; \cdots | \PP_2  \]
we can use classical tricks such as replacing disequalities by equalities with 
an extra variable to obtain
\[ \cdots;  P_2 + (tB) (P_1-P_2)  ; \cdots | \PP_2; B(Bt-1); t(Bt-1)  \]
Note that this step may increase the length exponentially, but the number of
variables grows only polynomially.

In conclusion, we can reduce {\em general} program equivalence to  
deciding $ \PP_1'| \PP_2'\approx \QQ_1' | \QQ_2' $,
where $ \PP_1',\PP_2',\QQ_1' $ and $ \QQ_2' $ are all arithmetic programs.
Let $ \ell $ be the input size of the original programs, namely, the
sum of the sizes of $ \PP_1,\PP_2, \QQ_1 $ and $ \QQ_2 $.    
The new arithmetic programs have length $ \exp(\ell) $ (the output size of the 
reduction). They have $ \exp(\ell) $ many
polynomials, but number of variables polynomial in $\ell$. The degree of each 
polynomial is at most $ \exp(\ell) $.
This reduction is a slightly improved version of the reduction \cite{BJK} used 
to derive their doubly exponential algorithm to solve the general universal
equivalence. Using our new algorithm for computing zeta functions of varieties, 
we can thus achieve a singly exponential time complexity. 

Let us now detail a key trick behind our improved zeta function algorithm. 

\section{Effective Kronecker theorem over finite fields}
\label{sec:kronecker} 
A classical theorem of Kronecker \cite{Kr} says that any affine algebraic set 
defined by a system of $m$ polynomials 
in $n$ variables over an algebraically closed field $K$ can be set theoretically defined by a system of 
$n+1$ polynomials in $n$ variables over the same field $K$.  Kronecker stated his theorem without a 
detailed proof; see \cite{Pe} for a self-contained proof. The theorem, as stated, is actually true for 
any infinite field $K$, not necessarily algebraically closed. But it fails for 
the finite field $\FF_q$, which is 
our main concern here.  In this section, we  follow the ideas in \cite{Pe} to show that Kronecker's theorem remains true for a 
finite field $\FF_q$ if $q$ is suitably large and we give an effective version of it tailored for our algorithmic application.   

Recall that if $I$ is an ideal in the commutative ring $\FF_q[x_1,..., x_n]$, then its {\em radical ideal} is defined as 
$\sqrt{I} =\left\{ \left. f \in \FF_q[x_1,..., x_n] \; \right| f^i \in I 
\ {\rm for ~some } \ i\geq 1\right\}$. 
It is then clear that the two ideals $I$ and $\sqrt{I}$ have the same set of 
$\FF_{q^k}$-rational points for every $k$.  
In particular, they have the same zeta function. 

\begin{theorem}[Affine version] 
\label{thm:kro} 
Let  $f_i \in \FF_q[x_1,\ldots, x_n]$ with $\deg(f_i)\leq d$ for all $1\leq i\leq m$. 
Assume that $q> (n+1)d^n$. Then there is a deterministic algorithm with 
running time 
$m(nd^n\log q)^{O(n)}$ which finds $n+1$ polynomials $g_j \in \FF_q[x_1,\ldots, x_n]$ 
with $\deg(g_j) \leq d$ for all $1\leq j\leq n+1$ such that their radical ideals are the same: $\sqrt{(f_1, \ldots,f_m)} = \sqrt{(g_1, \ldots, g_{n+1})}$. 
\end{theorem}

This theorem follows immediately upon dehomogenizing the following 
homogeneous version. 
\begin{theorem}\label{homver}
  [Homogeneous version] Let  $f_i \in \FF_q[x_1,\ldots, x_n]$ be homegenous polynomials of degree $d$ for all $1\leq i\leq m$. 
Assume that $q> nd^{n-1}$. There is a deterministic algorithm with running time 
$m(nd^n\log q)^{O(n)}$ which finds $n$ homogenous polynomials $g_j \in \FF_q[x_1,\ldots, x_n]$ 
of degree $d$ for all $1\leq j\leq n$ such that their radical ideals are the 
same: $\sqrt{(f_1, \ldots,f_m)} = \sqrt{(g_1, \ldots, g_n)}$. 
\end{theorem}

{\sl Proof of Theorem~\ref{homver}}. If $m\leq n$, the theorem is trivial as we can just take $g_j=f_j$ for $j\leq m$ and $g_j=f_1$ for $j>m$. 
We now assume that $m> n$. By induction, it is enough to prove the case $m=n+1$. Now, the $n+1$ 
polynomials $\{f_1,..., f_{n+1}\}$ in $n$ variables are algebraically dependent over $\FF_q$. That is, 
there is a non-zero homogenous polynomial $A_M(y_1,\ldots, y_{n+1})$ of some positive degree $M$ 
in $\FF_q[y_1,\ldots, y_{n+1}]$ such that 
$$A_M(f_1,\ldots, f_{n+1})=\sum_{k_1+\cdots +k_{n+1}=M} A_{k_1,\ldots, k_{n+1}}f_1^{k_1}\cdots f_{n+1}^{k_{n+1}} = 0.$$
This polynomial relation gives a homogenous linear system over $\FF_q$ with ${M+n\choose n}$ variables $A_{k_1,\ldots, k_{n+1}}$ 
and ${Md+n-1\choose n-1}$ equations. 
If ${M+n\choose n} >  {Md+n-1\choose n-1}$, 
the homogenous linear system will have a non-trivial solution. 
Now, choose $M= nd^{n-1}$. It is clear that $Md+i\leq d(M+i)$ for all $i\geq 0$ and 
$$\frac{{M+n\choose n}}{{Md+n-1\choose n-1}}= \frac{M+n}{n} \prod_{i=1}^{n-1}\frac{M+i}{Md+i} \geq \frac{M+n}{n}\left(\frac{1}{d}\right)^{n-1}>1.$$
Solving the linear system which takes time at most 
$$\left({M+n\choose n}\log q\right)^{\omega} = (M\log q)^{O(n)} = 
(nd^n\log q)^{O(n)},$$ 
(with $\omega<2.373$ the matrix multiplication exponent \cite{williams}), we 
can then clearly find a non-trivial solution 
$A_{k_1,\ldots, k_{n+1}}\in \FF_q$, with $k_1+\cdots +k_{n+1}=M$.

Next, we would like to make an invertible $\FF_q$-linear transformation 
$$f_u = \sum_{v=1}^{n+1} b_{u,v} g_v, \ b_{u,v}\in \FF_q, \ u=1,2,\ldots, n+1$$
such that when $A_M(f_1,\ldots, f_{n+1})$ is expanded as a polynomial in 
$\{g_1\cdots, g_{n+1}\}$ under the above linear tranformation, 
the coefficient of $g_{n+1}^M$ is non-zero. Such an invertible linear transformation may not 
exist if $q$ is small. We shall prove that it does exist if $q> M=nd^{n-1}$:  
Expand and write 
$$A_M(f_1,\ldots, f_{n+1}) = \sum_{k_1+\cdots +k_{n+1}=M} B_{k_1,\ldots, k_{n+1}}g_1^{k_1}\cdots g_{n+1}^{k_{n+1}}.$$ 
One checks that the coefficient of $g_{n+1}^M$ is 
$$\sum_{k_1+\cdots +k_{n+1}=M} 
A_{k_1,\ldots, k_{n+1}}b_{1. n+1}^{k_1}\cdots b_{n+1,n+1}^{k_{n+1}}= A_M(b_{1,n+1},\ldots, b_{n+1, n+1}).$$
This is a non-zero homogeneous polynomial in the $(n+1)$ variables $b_{u, n+1}$ ($1\leq u\leq n+1$) 
of degree $M$ with coefficients in $\FF_q$. Since $M<q$, the non-zero polynomial $A_M(y_1,\ldots, y_{n+1})$ 
is not the zero function on $\FF_q^{n+1}$. Now, a non-zero univariate polynomial $h(x)$ over $\FF_q$ of degree at most $M$ 
has at most $M$ roots in $\FF_q$. By trying at most $M+1\leq q$ elements of $\FF_q$, we find a non-root of $h(x)$ in $\FF_q$. 
Recursively applying this observation to the non-zero leading coefficient  (with respect to any one variable) of the non-zero polynomial $A_M(y_1,\ldots, y_{n+1})$, 
we find a non-zero vector  $(b_{1, n+1}, \ldots, b_{n+1, n+1})\in \FF_q^{n+1}$ such that 
$$c:= A_M(b_{1, n+1},\ldots, b_{n+1, n+1}) \in \FF_q^*.$$
This takes at most $(M+1)^{n+1} = (nd^{n-1}+1)^{n+1}$ trials. 
The non-zero vector $(b_{1, n+1}, \ldots, b_{n+1, n+1})$ can be easily extended to an invertible 
square matrix $(b_{u, v}) \in {\rm GL}_{n+1}(\FF_q)$. For instance, if $b_{n+1, n+1}\not=0$, then we can simply take 
$b_{u,v} = 0$ for $u \not=v$ and $1\leq v\leq n$, and $b_{u,v}=1$ if $u=v\leq n$.  In this way, we obtain the desired 
invertible transformation. 

Now, write our established polynomial relation in the form  
$$A_M(f_1,\ldots, f_{n+1}) = cg_{n+1}^M + G_1(g_1,\ldots, g_n) g_{n+1}^{M-1} +\cdots +G_{M}(g_1,\ldots, g_n)=0,$$
where $G_i(g_1,\ldots, g_n)$ is a homogenous polynomial in $\{g_1,\ldots, g_n\}$ of degree $i$ for $1\leq i\leq M$.  
Since the leading coefficient $c$ is not zero, we deduce that 
$g_{n+1}^M \in (g_1,\ldots, g_n)$. 
It follows that 
$$\sqrt{(f_1, \ldots,f_{n+1})} = \sqrt{(g_1,\ldots, g_{n+1})}= \sqrt{(g_1,\ldots, g_{n})}.$$
The theorem is proved.  \qed

\section{The Computation of Zeta Functions: Proving Theorem \ref{thm:main}}
Let $F$ be the following polynomial system with $m$ equations 
and $n$ variables over $\FF_q$: 
$$F(x_1,\ldots, x_n)=(f_1(x_1,\ldots, x_n),\ldots, f_m(x_1,..., x_n)),$$
where each $f_i \in \FF_q[x_1,\ldots, x_n]$ is a polynomial in $n$ variables of degree at most $d$. 
To compute the zeta function $Z(F, T)$, we need the following explicit degree 
bound of Bombieri: 
\begin{lemma} \cite{Bo} The total degree (the sum of the degrees for numerator 
and denominator) of the zeta function $Z(F, T)$ is bounded by $(4d+5)^{2n+1}$. 
\qed 
\end{lemma} 

Note that this total degree bound is independent of $m$. This already suggests 
the possibility of improving the dependence on $m$ in earlier algorithms for 
computing zeta functions. By applying our effective Kronecker theorem 
(Theorem \ref{thm:kro}), we are now ready to prove our main result. 

\medskip 
\noindent 
{\bf Proof of Theorem \ref{thm:main}:}  
If $q> (n+1)d^{n}$ then we can apply the affine effective Kronecker theorem in 
the previous section to replace the 
large polynomial system $F$ by a smaller polynomial system 
$G= (g_1(x_1,\ldots, x_n),\ldots, g_{n+1}(x_1,..., x_n))$, 
where each $g_j \in \FF_q[x_1,\ldots, x_n]$ is a polynomial in $n$ variables of degree at most $d$. 
The smaller system $G$ can be constructed in time 
$$m(nd^n\log q)^{O(n)} = m(n^nd^n\log q)^{O(n)},$$ 
thanks to Theorem \ref{thm:kro}. 
The two systems $F$ and $G$ have the same number of solutions over every extension field 
$\FF_{q^k}$. In particular, their zeta functions are the same, namely, 
$Z(F, T) = Z(G, T)$. 
Now, by the algorithms in \cite{Ha}, the zeta function $Z(G, T)$ can be computed in time 
$$2^{n+1} p ((n+1)^nd^n\log q)^{O(n)} = p (n^nd^n\log q)^{O(n)}.$$
Thus, the zeta function $Z(F,T)$ can be computed in time 
$$mp (n^nd^n\log q)^{O(n)}.$$

If $q\leq  (n+1)d^{n}$, we cannot apply the effective Kronecker theorem 
directly. So we use a somewhat different argument instead. 
Let $B=(4d+5)^{2n+1}$ be the upper bound in Bombieri's lemma. By \cite{Wa}, it is enough to 
compute the following $B$ numbers 
$$N_{q^k}(F), \ k=1, 2, ..., B.$$
If $q^k \leq (n+1)d^{n}$, namely, $k \leq  \log ((n+1)d^{n})/\log q$, we use the trivial exhaustive search algorithm 
to compute $N_{q^k}(F) $. For each such $k$, this takes time 
$$q^{k(n+1)} m (d^n \log q)^{O(1)} \leq 
((n+1)d^n)^{n+1}m (d^n \log q)^{O(1)}= m(n^nd^n\log q)^{O(n)}.$$
If $q^k \geq (n+1)d^{n}$, namely, $\log ((n+1)d^{n})/\log q \leq k\leq B$, then 
we can apply the effective Kronecker theorem 
to the system over the extension field $\FF_{q^k}$ 
to produce a new system
$$G_k= (g_{k,1}(x_1,\ldots, x_n),\ldots, g_{k,n+1}(x_1,..., x_n)),$$
where each $g_{k,j} \in \FF_{q^k}[x_1,\ldots, x_n]$ is a polynomial in $n$ variables of degree at most $d$. 
Now, 
$$N_{q^k}(F) = N_{q^k}(G_k).$$
The system has only $n+1$ equations and thus the number $N_{q^k}(G_k)$ (in fact the full zeta function 
of $G_k$ over $\FF_{q^k}$) can be computed 
by \cite{Ha} in time 
$$2^{n+1}p(k(n+1)^nd^n\log q)^{O(n)}= p(Bn^nd^n\log q)^{O(n)} = p(n^nd^n\log q)^{O(n)}.$$
Thus, the total time to compute $Z(F, T)$ is bounded by 
$$B mp (n^nd^n\log q)^{O(n)} = mp(n^nd^n\log q)^{O(n)}. \text{ \qed} $$

\end{document}